\documentclass{amsproc}

\usepackage{amssymb,xypic}
\usepackage[dvips]{graphicx}

\newtheorem{theorem}{Theorem}[section]

\newtheorem{prop}[theorem]{Proposition}

\theoremstyle{definition}
\newtheorem{definition}[theorem]{Definition}
\newtheorem{eg}[theorem]{Example}
\theoremstyle{remark}

\numberwithin{equation}{section}

\let\To=\longrightarrow
\let\into=\hookrightarrow
\newcommand\Pee{\mathbb P}
\newcommand\Ltimes{\mathbin{\buildrel{\mathbf L}\over\otimes}}
\newcommand\Lf{\mathbf Lf^*}
\newcommand\OO{\mathcal O}
\newcommand\C{\mathbb C}
\newcommand\F{\mathcal F}
\newcommand\G{\mathcal G}
\newcommand\Ll{\mathcal L}
\newcommand\RR{\mathbb R}
\newcommand\R{\mathbf R\hspace{1pt}}
\newcommand\res{\arrowvert_}

\newcommand\ve{^\vee}
\newcommand\Hom{\mathrm{Hom\,}}
\newcommand\Ext{\mathrm{Ext}}
\newcommand\id{\mathrm{id\,}}

\begin{document}

\title{Mirror symmetry and actions of braid groups on derived categories}

\author{R.\,P. Thomas}
\address{Mathematical Institute, 24--29 St Giles', Oxford OX1 3LB. UK}
\email{thomas@maths.ox.ac.uk}
\thanks{The author is supported by Hertford College, Oxford}

\date{30 April 99}

\begin{abstract}
After outlining the conjectural relationship between the
conjectural mirror symmetry programmes of Kontsevich and
Strominger-Yau-Zaslow, I will describe some natural consequences of this
which are proved from scratch in joint work with Mikhail Khovanov and
Paul Seidel.
Namely, actions of braid groups are found on derived categories of
coherent sheaves, dual to Seidel's braid group of symplectic
automorphisms (``generalised Dehn twists'') of the mirror. These
Fourier-Mukai transforms, one for every ``spherical'' element of the
derived category (a simple rigid sheaf on a Calabi-Yau 3-fold is an
example) are closely related to the ``mutations'' of exceptional
bundles on Fanos. Examples of conjecturally mirror dual group
actions on triangulated categories are drawn from smoothings
and resolutions of 2 and 3 dimensional singularities.
\end{abstract}

\maketitle

\section{Introduction}

In this talk I will attempt to do two things: review a conjectural
picture of mirror symmetry, which will no doubt crop up in many other
talks, and explain some of its consequences which are proved (independently
of the conjectures) in joint work with Mikhail Khovanov and Paul Seidel
\cite{ST,KS}. This is not the mirror symmetry of Gromov-Witten
invariants and variations of Hodge structure (yet) but the more
fundamental programmes of Kontsevich \cite{K} and Strominger-Yau-Zaslow
\cite{SYZ} which should eventually lead back to the traditional
predictions of mirror symmetry.

Section \ref{KSYZ} explains the two programmes and their supposed link
in the language of Fourier-Mukai transforms,
which are reviewed and explained in Section \ref{FM}. Though many of
the objects used in the models are yet to be defined, we can
attempt to deduce some consequences and prove these. In
particular, symplectomorphisms of a Calabi-Yau manifold should induce
autoequivalences of the bounded derived category $D^b(X)$ of coherent
sheaves on the mirror $X$. We find many such autoequivalences, and
actions of the braid group on $D^b(X)$ in particular cases (such as
on resolutions of singularities). These should be mirror dual to braid
groups of symplectomorphisms that arise from configurations of
Lagrangian spheres. These spheres are often the vanishing cycles of the
smoothing of another, dual, singularity, and we discuss this mirror symmetry
of singularities. In Section \ref{mut} we mention
briefly the relationship of this work to mutations of bundles on Fano
manifolds, and finally in Section \ref{fid} we outline the
faithfulness of the braid group action.

\textbf{Acknowledgements.} My main debts are to Mikhail Khovanov and
Paul Seidel who allowed me to talk about this work as if it
were my own. Kontsevich and Bridgeland and Maciocia have also discussed
the twist (\ref{twist}). There are many more acknowledgements in
\cite{ST}. I would also like to thank the organisers
of the 1999 Harvard Winter School on Mirror Symmetry for inviting me
to take part.

\section{Fourier-Mukai transforms} \label{FM}
\subsection*{Function transforms}

We begin by introducing Fourier-Mukai transforms by analogy with
function transforms. Suppose we are given a family $F_p$ of
(complex-valued, say) functions or distributions on some space $V$,
parametrised by some dual space $\widehat V\ni p$.

Let $F:\,V\times\widehat V\to\C$ be the universal function (with
$F\res{V\times\{p\}}=F_p$).
For instance we could take $V=\RR^n,\ \widehat V=(\RR^n)^*$ with
$F(x,p)=e^{ip.x}=F_p(x)$. Similarly on any manifold
$V$ we can take $\widehat V=V$ and $F=\delta_\Delta$, the Dirac-delta
of the diagonal $\Delta$, with $F_p(x)=\delta_p(x)$.

Suppose now that $\{F_p\,:\,p\in\widehat V\}$ span (some class of)
functions $V\to\C$, and are orthonormal with respect to some inner
product (e.g. $\int F_p\overline F_q=\delta_{pq}$ in a distributional
sense). Then for any $f$ in this class, $f$ is built up from the
$F_p$\,s, with the coefficient of $F_p$ being $\widehat
f(p)=\int_Vf\overline F_p$, and we sum $\widehat f(p)F_p$ over
$p$ to regain $f$. That is,
$$
f=\int_{\widehat V}\left[\int_Vf\overline F_p\right]F_p\ dp,
$$
so $f=(\widehat f\,)\ve$, where $\ve$ is the dual transform
$g\ve(x)=\int_{\widehat V}g(p)F(p,x)\,dp$.

For instance for $F=e^{ip.x}$ we get the Fourier Inversion theorem
(though theorem is a little strong for the above treatment) for the
Fourier transform $\ \widehat{}\ $ and its inverse $\ve$. Similarly taking
$F=\delta_\Delta$ on $V\times V$ gives two applications of the
identity map: $\widehat f(p)=\int_Vf\delta_p=f(p)$.

The way to look at this most relevant to sheaves is via the diagram
$$
\spreaddiagramcolumns{-1.8pc}
\spreaddiagramrows{-1pc}
\diagram
& V\times\widehat V \dlto_(.6){\pi_1}\drto^(.6){\pi_2} \\ V &&
\widehat V.
\enddiagram
$$
Then the transforms are
\begin{eqnarray} \nonumber
\widehat f\!&=&\!(\pi_2)_*\big[\pi_1^*f\,.\,\overline F\,\big], \\
g\ve\!\!\!&=&\!(\pi_1)_*\big[\pi_2^*g\,.\,F\,\big].  \label{ft}
\end{eqnarray}
Thus we pull up a function from $V$ to the product, multiply by
$\overline F$, and push down the fibres $V_p$ of $\pi_2$
($\pi_*$ means integrate
down the fibres of $\pi$): at each $p$ we integrate $f$ against
$\overline F_p$ to get $\widehat f(p)$. Notice also that we can
recover the functions $F_p$ parametrised by $p$ as the inverse
transform of the Dirac-delta $\delta_p$,
\begin{equation} \label{p}
F_p=\delta_p\ve.
\end{equation}

We can also do this in a family. For instance on $V\times Z$ we can
carry an extra parameter $z\in Z$ and take the usual transform at each
point $z$, i.e. the transform just with respect to the variables $x$
and $p$,
\begin{equation} \label{rel}
f(x,z)\stackrel{\widehat{}}{\mapsto}\widehat f(p,z).
\end{equation}
This is just a relative transform on $(V\times Z)\times_Z(\widehat
V\times Z)=V\times\widehat V\times Z$. Alternatively we can write it
as a strict transform on $(V\times Z)\times(\widehat V\times Z)$. The
function parametrised by $(p,z)\in\widehat V\times Z$ is
$F_p.\delta_z$, where $F_p$ is (the pull-back to $V\times Z$
of) one of the family of functions on $V$, and $\delta_z$ is
the Dirac-delta concentrated on $V\times\{z\}$. Again these functions
$$
\{F_{p,z}=F_p.\delta_z\,:\,(p,z)\in\widehat V\times Z\}
$$
form an orthonormal set, giving a transform with universal function
$F\delta_\Delta$, where $F$ is the pull-back of the universal function
on $V\times\widehat V$ and $\Delta=V\times\widehat V\times
Z\subset(V\times Z)\times(\widehat V\times Z)$ is the diagonal. The
transform $\widehat f(p,z)$ is of course the same as (\ref{rel}).

Similarly we can have a non-trivial family $X\to Z$, with dual family
$\widehat X\to Z$, and a relative transform on $X\times_Z\widehat X$
or a strict transform on $X\times\widehat X$. Again they are the same,
with the universal function on $X\times_Z\widehat X$ pushed forward
via the diagonal map $X\times_Z\widehat X\into X\times\widehat X$ to a
function supported on this diagonal. This will seem much more natural
in the setting of sheaves.

\subsection*{Sheaf transforms}

Suppose now we have a family $F_p$ of vector bundles, or coherent
sheaves, or complexes of sheaves, on some complex manifold $X$,
parametrised by some dual variety $\widehat X\ni p$. We will also
need to assume the existence of a universal object $F$ on $X\times
\widehat X$ whose restriction to each $X\times\{p\}$ is $F_p$.

We will define the inverse transform first (denoted here by the
unfortunate notation $\ve$; this is not the dual of a sheaf but its
transform). We would like to think of a sheaf/bundle/complex $\F$ on 
$\widehat X$ as giving a distribution of the $F_p$-components of a
sheaf $\F\ve$ on $X$, i.e. $\F\ve$ should be the sum over $p\in\widehat
X$ of the sheaves $\F\res p\otimes F_p$. So with respect to the
diagram
$$
\spreaddiagramcolumns{-1.8pc}
\spreaddiagramrows{-1pc}
\diagram
& X\times\widehat X \dlto_(.6){\pi_1}\drto^(.6){\pi_2} \\ X &&
\widehat X
\enddiagram
$$
we would like to set, as in (\ref{ft}),
$$
\F\ve=(\pi_1)_*\big[(\pi_2^*\F)\otimes F\big].
$$
Since a sheaf assigns groups of sections to open sets, the right
notion of pushdown $(\pi_1)_*$ (the sum of all of the sections over
$p\in\widehat X$, the analogue of integration down the fibres) is to
assign to an open set $U$ all sections on
$\pi_1^{-1}(U)=U\times\widehat X$. In fact numbers of sections down
fibres can jump, or be zero, due to the presence of higher
cohomology. So we actually take the right derived
functor $\R(\pi_1)_*$ described in my last talk -- this is the
complex of sheaves obtained by resolving by injective sheaves and
applying $(\pi_1)_*$ to this complex. For similar reasons the tensor
product should be the derived tensor product $\Ltimes$, and
the restriction to fibres $\otimes\OO_{X\times\{p\}}$ we have
mentioned before should also be taken in this derived sense.

Thus the correct setting is the bounded derived category of coherent
sheaves $D^b$, with a universal object $F\in D^b(X\times\widehat
X)$. We set, for $\F\in D^b(\widehat X)$,
\begin{equation} \label{F}
\F\ve=\R\pi_{1\,*}\big[(\pi_2^*\F)\Ltimes F\big];
\end{equation}
in most of the simple cases we shall be interested in one can think of
vector bundles rather than complexes of sheaves, and normal pushdown
and tensor product.

For instance the inverse transform of the structure sheaf $\OO_p$ of a point
$p\in\widehat X$ is
\begin{equation} \label{op}
\OO_p\ve=\R\pi_{1\,*}\big[\OO_{X\times\{p\}}\Ltimes
F\big]=\pi_{1\,*}F_p=F_p,
\end{equation}
the (derived) restriction of $F$ to $X\times\{p\}$, i.e. the object of
$D^b(X)$ parametrised by $p\in\widehat X$; compare
(\ref{p}). Similarly the inverse transform of $\bigoplus_i\OO_{p_i}$
is $\bigoplus_iF_{p_i}$.

The transform from $D^b(X)$ to $D^b(\widehat X)$ is given by the
formula (c.f. (\ref{ft}))
$$
\widehat{}\ =\R\pi_{2\,*}\big[(\pi_1^*(\,\cdot\,)\Ltimes F^*\big]\,[\,n\,],
$$
where here we have denoted by $F^*$ the derived dual of $F$,
$\R\mathcal Hom(F,\OO_{X\times\widehat X})$, and $n=$\,dim\,$X$.
Thus for $\F\in D^b(X)$ we have
$$
\widehat\F=\R\pi_{2\,*}\big[\R\mathcal Hom(F,\pi_1^*\F)\big]\,[\,n\,].
$$
For $\ve$ and $\ \widehat{}\ $ to be actual inverses (equivalences is the
correct categorical notion) we need the objects
$F_p=F\Ltimes\OO_{X\times\{p\}}$ to be orthogonal \cite{BO},
\begin{equation} \label{1}
\R\Hom(F_p,F_q)\cong0, \quad p\ne q,
\end{equation}
orthonormal, in the sense that they are \emph{simple},
\begin{equation} \label{2}
\Hom(F_p,F_p)=\C\,.\,\mathrm{id},
\end{equation}
and to satisfy the dimension and partial Calabi-Yau conditions \cite{Ma}
\begin{eqnarray}
\mathrm{dim}\ X=\mathrm{\,dim}\ \widehat X, \nonumber \\
F_p\otimes\omega_X\cong F_p \quad \forall p. \label{3}
\end{eqnarray}
Recently Bridgeland \cite{Br} has shown that for $X$ and $\widehat X$
smooth projective varieties these conditions (\ref{1}), (\ref{2}),
(\ref{3}) are also sufficient for the
Fourier-Mukai transform given by $F\in D^b(X\times\widehat X)$ to be an
equivalence. Here $\omega_X$ is the canonical bundle of $X$, and we
are requiring this to be trivial on the support of each $F_p$; recall
that $X$ is Calabi-Yau if $\omega_X$ is globally trivial.

\begin{eg}
Take $\widehat X=X$ and $F=\OO_\Delta$, the structure sheaf of the
diagonal parametrising the sheaves $F_x=\OO_x$. Then the transform is
the identity, and the partial Calabi-Yau condition (\ref{3}) is
vacuous since $\OO_x\otimes\omega_X$ is trivially isomorphic to
$\OO_x$.
\end{eg}

\begin{eg}
The original Fourier-Mukai transform \cite{Mu}. Let $T$ be an elliptic curve,
and let $\widehat T\ni t$ be its Jacobian parametrising degree 0 line
bundles $L_t$ on $T$. We take $F=L$ to be a Poincar\'e line bundle on
$T\times\widehat T$ (for instance if we fix a basepoint $t_0\in T$
then we can identify $T$ with $\widehat T$, and $L$ is given by the
divisor $\Delta-(\{t_0\}\times T)\subset T\times T$).

Then the $L_t$\,s form an orthonormal set:
$$
H^0(L_s^*\otimes L_t)=H^1(L_s^*\otimes L_t)=\left\{
\begin{array}{ll}\C & s=t, \\ 0 & s\ne t. \end{array} \right.
$$
So as $T$ is Calabi-Yau we get an invertible transform
$$
\F\mapsto\widehat\F=\R\pi_{1\,*}\big[L\otimes\pi_2^*\F\big]
$$
with fibre over the generic point $t\in\widehat T$ (at which
base-change holds) ``$H^0(\F\otimes L_t)-H^1(\F\otimes L_t)$''
(really a complex with these two cohomologies). Thus, when $H^1$
vanishes, we replace $\F$ by its spectrum of sections $\{H^0(\F\otimes
L_t):t\in\widehat T\}$ of different twists. The inverse transform
uses the dual line bundle, and the set-up is symmetric between $T$ and
$\widehat T$.

Then using (\ref{op}) we see that the transform of the structure sheaf
$\OO_t$ of a point $t\in\widehat T$ is the corresponding line bundle
$L_t$, while $r$ points transform to the appropriate rank $r$ bundle
(which is the sum of $r$ line bundles). We depict this as follows,
drawing a basis of the fibre of the vector bundle over the right hand
torus. 
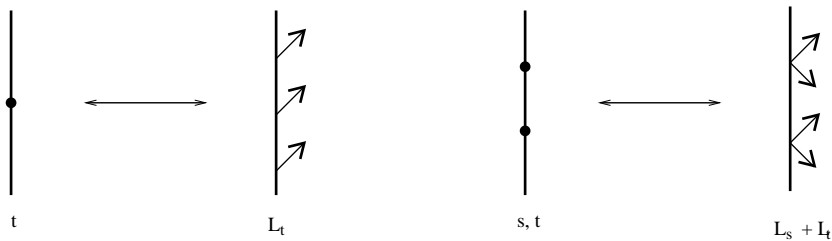
\vspace{-3mm}
\begin{figure}[h]
\begin{center}
\includegraphics[angle=0,width=11cm]{fm1.eps} \vspace{-4mm}
\end{center}
\end{figure}

Thus, inversely, we have an algebraic gadget that takes a degree 0
semistable bundle on $T$ and splits it into its constituent pieces
according to Atiyah's classification; for every $L_t$ factor we get a
point $t\in\widehat T$. The non-trivial extensions correspond to
structure sheaves of double points, etc.
\end{eg}

Similarly we can do this in a family. Consider an elliptic fibration
$X\to Z$, possibly with singular fibres. Suppose we also have a
section $s$ so that we can identify $X$ with its relative
Jacobian. Then there is a Poincar\'e sheaf $\Ll$ on $X\times_Z X$
corresponding to the Weil divisor
$\Delta-s(Z)\times_ZX$. Rather than doing a relative transform
on this singular space, we do a strict transform on $X\times X$ by
pushing forward $\Ll$ to $\iota_*\Ll$ via $\iota:\,X\times_Z X\to X\times
X$. Then setting
$$
\widehat\F=\R\pi_{1\,*}\big[\iota_*\Ll\otimes\pi_2^*\F\big]
$$
we get a transform $D^b(X)\to D^b(X)$, a family version of the
previous example, giving pictures like:

\begin{figure}[h] \nopagebreak
\begin{center}
\includegraphics[angle=0,width=11cm]{fm2.eps} \vspace{-2mm}
\end{center}
\end{figure}
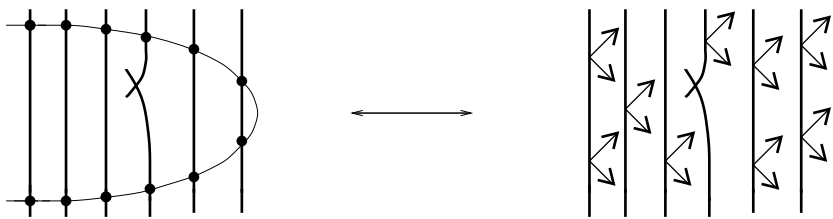
\hspace{-5mm}Multisection $C$ (+ line bundle on it) $\quad
\longleftrightarrow\quad$\ \ Vector bundle, deg 0 on fibres

\noindent where the covering degree $r$  of the multisection (its
intersection with the fibre) is the rank $r$ of the vector bundle.

Thus vector bundles correspond to ``spectral covers'' in the usual way
-- the cover gives $r$ points on each fibre yielding the sum of $r$
line bundles on the dual fibre. Globally we patch these together in
the base direction using the line bundle on $C$, and get a non trivial
vector bundle which is not globally a sum of line bundles (because $C$
is not globally $r$ disjoint sections).

The rigorous statement is that there is an autoequivalence of $D^b(X)$
given by this transform. It is invertible since the universal sheaf
$\iota_*\Ll$ is concentrated on torus fibres, about which $X$ is
locally Calabi-Yau so that $\iota_*\Ll\otimes\omega_X\cong\iota_*\Ll$.
For us, the point is simply that there is an algebraic gadget that converts
spectral covers like $C$, with line bundles on them, into vector bundles
(degree 0 and semistable on the fibres), and vice-versa. It does this
for free, without any analysis of the singular fibres, since $X\times X$ is
smooth (even though $X\times_Z X$ is not). This is at the expense of
introducing the derived category, of course, but as we see here in
many cases the transform takes sheaves to sheaves instead of a
complex of sheaves.

\section{Kontsevich vs. Strominger-Yau-Zaslow} \label{KSYZ}

We now turn to mirror symmetry and the two competing conjectural
theories. Strominger-Yau-Zaslow \cite{SYZ} suggest that (in an
appropriate complex structure) a Calabi-Yau $n$-fold $X$ should admit
a fibration by special Lagrangian tori ($T^n\into X$ is Lagrangian if
the symplectic (K\"ahler) form restricts identically to zero on $T^n$,
and special if the restriction of the imaginary part of the
Calabi-Yau $(n,0)$-form is
identically zero) with a special Lagrangian section. In this
case the mirror $\widehat X$ should be the dual torus fibration.

Kontsevich \cite{K}, on the other hand, conjectures that there should
be a natural exact equivalence of triangulated categories (an exact
equivalence is one which preserves the distinguished triangles) between
$D^b(X)$ and $D^b(\mathrm{Fuk}\,(\widehat X))$. The second category here
is the derived category of the Fukaya category of the symplectic manifold
$\widehat X$, and as such is not yet properly constructed (though
see Fukaya's talk for more progress on this). Other talks will explain
more about the Fukaya category; all we need to know is that it is
constructed from only the symplectic geometry of $\widehat X$, using
(graded) Lagrangian submanifolds, local systems on them, and their
Fukaya-Floer homology. In particular every Lagrangian cycle $L$ in
$\widehat X$ (plus a grading), with a flat line bundle on it, should
give an object in $D^b(\mathrm{Fuk}\,(\widehat X))$. The corresponding
object $\F_L$ in $D^b(X)$ should have the same Homs, so that $HF^*(L,L)$
should be quasi-isomorphic to $\R\Hom(\F_L,\F_L)$.

The (again conjectural) correspondence between the two pictures is now
folklore and has been discussed by many people, see for instance
\cite{AP,Ty}. The basic idea is that a Lagrangian multisection $L$ in the
fibration $\widehat X$ (with a flat line bundle on it, and intersection with
the fibre $r$) should correspond to a rank $r$ holomorphic vector
bundle on $X$ by an analytical version of the Fourier-Mukai transform,
giving a diagram like the one above. That is, the
intersection of $L$ with a fibre gives $r$ points corresponding to $r$
line bundles on the dual torus, as before. \emph{Special} Lagrangian
sections should perhaps correspond to bundles with
Hermitian-Yang-Mills connections (i.e. \emph{stable} bundles) as
suggested in \cite{Va} (both special and stable are
stability conditions on the objects on either side under a natural
group action, so this makes sense).

In 2 complex dimensions we can rigorously carry out this procedure,
since we have the tools of algebraic geometry and Fourier-Mukai to
deal with the singularities. We also have Yau's solution of the Calabi
conjecture, which gives a hyperk\"ahler metric on a compact Calabi-Yau
surface (i.e. a $K3$ or $T^4$). Thus for every complex structure $I$
there is a conjugate quaternionic partner $J$, and if we rotate the
complex structure from $I$ to $J$ then the special Lagrangian cycles
become the complex curves on $X$. Thus, after this hyperk\"ahler
rotation, the SYZ conjecture is concerned with an elliptically fibred
surface with a section, and the mirror should be the dual
fibration. Thus in this dimension the mirror is topologically the same
as $X$ (in 3 dimensions there is a topology change at the singular fibres)
and the correspondence we want is precisely the Fourier-Mukai
transform described before. This gives an equivalence of categories
$D^b(X)\to D^b(X)$, which is the appropriate even dimensional mirror
symmetry, taking 0, 2, 4 branes (points, 2-cycles, and vector
bundles) to even dimensional 0, 2, 4 branes (3 dimensional mirror
symmetry as we have described should take 0, 2, 4, 6 branes to 3 branes).

We can ask what else we get from the Fourier-Mukai transform. By
(\ref{op}) a fibre,
or more generally the pushforward of a line bundle on a fibre, gets
taken to the corresponding point in the dual fibre:

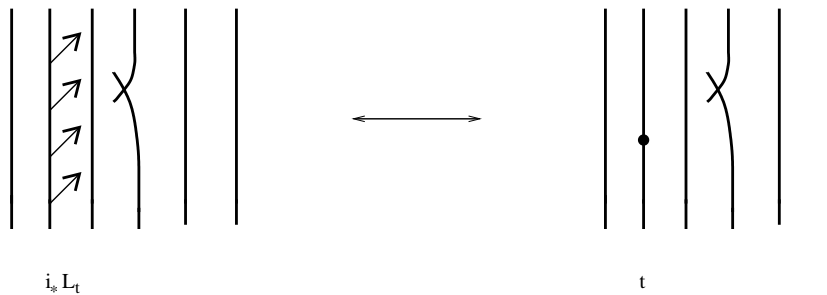
\begin{figure}[h]
\begin{center}
\includegraphics[angle=0,width=11cm]{fm3.eps}
\end{center}
\end{figure}
This was of course the original idea of SYZ: that the moduli of
special Lagrangian tori, plus flat line bundles, on $\widehat X$
should be isomorphic to the moduli of points on $X$, i.e. to
$X$ itself.

On non-hyperk\"ahler manifolds we cannot deduce anything much
about either programme without hard analysis, and that is best left to
other speakers here. But we can deduce some consequences, and try to
prove these. The surface
case which could be dealt with rigorously showed one thing: while the
mirror $\widehat X$ was isomorphic to the original $X$, the mirror map
was certainly \emph{not} the identity: in fact it took points to
fibres plus line bundles, rather than points, and induced a non
trivial map on $H^*(X)$ on taking Mukai vectors.

In particular we see that \emph{mirror symmetry is not
functorial on points} (a phrase I learnt from Paul Seidel); in fact, as
Kontsevich envisaged, \emph{(graded) symplectic automorphisms of the
mirror $\widehat X$ should not induce holomorphic automorphisms of
$X$, but autoequivalences of its derived category $D^b(X)$}.

This is what we shall concentrate on.

\section{Autoequivalences of derived categories and braid groups} \label{meat}

\subsection*{Autoequivalences of derived categories}

What are the autoequivalences of $D^b(X)$\,? There are the obvious ones
given by translation of complexes $[\,n\,]$, tensoring with line
bundles $\otimes L$, and those pulled back from automorphisms of
$X$. Bondal and Orlov \cite{BO} have shown that this is all of them
for $X$ smooth and projective with ample canonical bundle or anticanonical
bundle. In fact for there to be any more $X$ must be partially Calabi-Yau:
Orlov \cite{Or} has shown that any autoequivalence is set up
by an object $F\in D^b(X\times\widehat X)$ on the product as a
Fourier-Mukai transform (\ref{F}), which by (\ref{3}) must then satisfy
$F_p\otimes\omega_X\cong F_p$ for all $p\in X$. And for Calabi-Yau
manifolds Kontsevich's proposal predicts there should be many such
Fourier-Mukai transforms; in particular there should be one for every
(graded) symplectomorphism of the mirror $\widehat X$.

For instance Seidel \cite{S} has shown that given any Lagrangian
sphere $L\cong S^n$ in a symplectic manifold we may construct a
symplectomorphism -- the generalised Dehn twist about $L$. This is
given as monodromy around a degeneration of the manifold in which the
sphere is collapsed to (becomes the vanishing cycle of) an ordinary
double point. Alternatively we can glue in the following local model
twist on $T^*S^n$. Give $T^*S^n$ its
standard symplectic structure and metric. Then $\mu(\xi)=|\xi|$, the
length of a cotangent vector $\xi$, gives a Hamiltonian function which
induces a circle action -- the flow from $\xi$ (considered as a
tangent vector using the metric) is the \emph{unit speed} geodesic
flow in the direction $\xi$ along $S^n$, lifted horizontally to
$TS^n\cong T^*S^n$ as $(x(t),\xi(t)=\dot x(t)|\xi(0)/\dot x(0)|)$.
This flow is clearly discontinuous across the zero section $S^n\subset
T^*S^n$ (it is unit speed in opposite directions as we pass through the
zero section) but since the geodesics have length $2\pi$, the flow
through time $\pi$ \emph{is}
continuous, and gives the antipodal map. So we may define the
generalised Dehn twist as the flow of any point $\xi$ through an angle
varying smoothly from $0=2\pi$ as $|\xi|\to\infty$ to $\pi$ at the
zero section $|\xi|=0$; see \cite{S}.

The Dehn twist has a canonical lift to the graded symplectomorphism
group of $\widehat X$, and so should be dual to a Fourier-Mukai
transform constructed from an element $\F_L$ of $D^b(X)$ (not
$D^b(X\times X)$, notice) dual to the Lagrangian $L$. Since Homs
should be the same on both sides, we know that $\R\Hom(\F_L,\F_L)$
should be isomorphic to $HF^*(L,L)\cong H^*(S^n)$. Thus we might expect to
be able to find an invertible Fourier-Mukai transform for every
\emph{spherical} $\F\in D^b(X)$, where

\begin{definition} \label{sph}
$\F\in D^b(X)$ is $N$-spherical if $\R\Hom(\F,\F)\cong H^*(S^N;\C)$, where
$N=\mathrm{\,dim\,}X$. That is
$$
\Ext^i(\F,\F)=\left\{\begin{array}{cl}\C & i=0,N, \\ 0 & i\ne0,N.
\end{array} \right.
$$
\end{definition}

Thus for instance a simple (End($\F)=\C\,.\,\id$), rigid
($\Ext^1(\F,\F)=0$) sheaf $\F$ on a smooth Calabi-Yau 3-fold is
3-spherical by Serre duality: $\Ext^i(\F,\G)\cong\,\Ext^{3-i}(\G,\F)^*$.

\begin{definition}
For an $N$-spherical $\F\in D^b(X)$ define the twist $T_\F\G$ of $\G\in
D^b(X)$ to be the cone (total complex) of the evaluation map
\begin{equation} \label{twist}
\F\Ltimes\R\Hom(\F,\G)\to\G.
\end{equation}
Here we should pick suitable resolutions so the above becomes a genuine
map of complexes, then take the total complex of this map. This defines
$T_\F\G$ only up to quasi-isomorphism (cones are not functorial in
$D^b(X)$); $T_\F$ can in fact be made into a functor $D^b(X)\to
D^b(X)$ \cite{ST}.
\end{definition}

So in simple cases $T_\F\G$ will be the kernel or cokernel of a map
$\F\otimes\Hom(\F,\G)\to\G$. For $X$ a smooth projective variety, it
is easy to see that $T_\F$ is the Fourier-Mukai transform given by the
object
\begin{equation} \label{obj}
\{\R\mathcal Hom(\pi_2^*\F,\pi_1^*\F)\to\OO_\Delta\}
\end{equation}
of $D^b(X\times X)$. Here $\pi_i$ is
projection to the $i$th factor $X$, the map is restriction to the
diagonal $\Delta$ followed by the trace map, and the braces mean we
take the cone (total complex) of the map. However, it is more
convenient for us to work with general twists in arbitrary derived
categories; thus for instance our results apply to non-compact schemes
without difficulty.

The categorical equivalent of the partial Calabi-Yau condition
$\F\otimes\omega_X\cong\F$ is the existence of a functorial duality
(Serre duality for us, or a local form of it in the non-compact case)
between $\R\Hom(\F,\G)$ and $\R\Hom(\G,\F)\,[N]$; equivalently the
pairing
\begin{equation} \label{dual}
\Ext^i(\F,\G)\otimes\Ext^{N-i}(\G,\F)\stackrel{\smallsmile}{\To}
\Ext^N(\F,\F)\to\C
\end{equation}
should be perfect, where the last map uses the fact that $\F$ is
$N$-spherical.

\begin{definition}
For $\F$ $N$-spherical (\ref{sph}) and with a duality (\ref{dual}), there
is a functor $T'_\F:\,D^b(X)\to D^b(X)$ \cite{ST} such that the
quasi-isomorphism class of $T'_\F\G$ is the total complex we denote
$$
\{\G\to\F\Ltimes\R\Hom(\F,\G)\,[N]\}
$$
given by dualising a map of chain complexes representing the evaluation
$$
\G\Ltimes\R\Hom(\G,\F)\to\F.
$$
\end{definition}

Again this is really a Fourier-Mukai transform with object the derived
dual of (\ref{obj}) shifted by $[N]$.

The last three definitions, and the theorems below, can of course be
formulated in the derived category of an arbitrary abelian category,
linear over a field $k$, having enough injectives, and  containing a
spherical $\F$ with a duality (\ref{dual}). We shall confine ourselves
to derived categories of coherent sheaves (this does not have enough
injectives and one must work with quasi-coherent complexes with
coherent cohomology in the usual way; see \cite{ST} for full details).

\begin{theorem}
For $\F$ spherical with a duality (\ref{dual}), $T_\F$ and $T'_\F$ are
inverses.
\end{theorem}

Again we are being sketchy here. The precise statement \cite{ST} takes
place in the derived category of a $k$-linear abelian category, and
``inverses'' means that $T_\F\circ T'_\F$ and $T'_\F\circ T_\F$ are both
naturally isomorphic to the identity functor.

Using Serre duality on a compact scheme, or a local form of it one can
prove for $\F$ compactly supported on a non-compact scheme, we get
an invertible Fourier-Mukai transform in each of the following examples.

\begin{itemize}
\item Any simple, rigid sheaf $\F$ on a Calabi-Yau 3-fold $X$. In particular
  the structure sheaf $\OO_X$ gives a canonical transform, called the
  reflection functor by Mukai, which should be mirror dual to the Dehn
  twist about the (conjectural) special Lagrangian $S^3$ zero section
  of the SYZ fibration of $\widehat X$.
\item Holomorphic $-2$-spheres $C$ in a complex surface give Fourier-Mukai
  transforms, taking $\F=\OO_C$. In particular we get Fourier-Mukai
  transforms for general-type surfaces containing such spheres, which
  might be surprising given the results of \cite{BO}. These surfaces
  are locally Calabi-Yau along the spheres. The induced action on
  cohomology is the Picard-Lefschetz reflection in the corresponding
  $-2$-vector.
\item Spheres with normal bundle
  $\OO_{\Pee^1}(-1)\oplus\OO_{\Pee^1}(-1)$ in 3-folds also give
  transforms in the same way.
\item Surfaces $S$ in 3-folds satisfying $h^{0,1}(S)=0=h^{0,2}(S)$,
  with the local Calabi-Yau condition $\nu_S\cong\omega_S$.
\end{itemize}

Seidel proved in \cite{S} that the Dehn twists along
$A_n$-chains of Lagrangian spheres satisfy the braid relations. That
is, if we have a chain of such spheres with consecutive pairwise
intersections one transverse point, then we get a homomorphism from the
braid group on $(n+1)$ strands $B_{n+1}$ into the
symplectomorphism group of the ambient manifold. Moreover he
showed the smoothing of an $A_n$ singularity
on a complex surface (such as the smoothing of the standard $SL(2,\C)$
quotient singularity $\C^{\,2}/\mathbb Z_n$) contains such a chain of
Lagrangian spheres.

For two Lagrangians $L_i$ intersecting transversely in a single point
we have $HF^*(L_1,L_2;\C)\cong\C$, so we define

\begin{definition}
A sequence of spherical objects $\F_i\in D^b(X),\ i=1,\ldots,n$ form
an $A_n$-chain if they satisfy
$$
\sum_k\mathrm{dim\,}\Ext^k(\F_i,\F_j)=\left\{\begin{array}{cl}
0 & |i-j|>1, \\ 1 & |i-j|=1. \end{array}\right.
$$
\end{definition}

Thus there are no Homs between distinct $F_i$\,s unless they are
consecutive, in which case there is a unique map (up to scale) in some
(arbitrary) degree.

\begin{theorem}
Given an $A_n$-chain of spherical objects $\F_i\in D^b(X)$ with
duality (\ref{dual}) there are the following natural isomorphisms
between the corresponding functors $T_i=T_{\F_i}$
$$ \begin{array}{ll}
T_i\,T'_i\cong\id\cong T'_i\,T_i \\
T_i\,T_j\,T_i\cong T_j\,T_i\,T_j\qquad & |i-j|=1, \\
T_i\,T_j\cong T_j\,T_i & |i-j|>1. \end{array}
$$
Thus they define a weak braid group action on $D^b(X)$. (We have not
checked if the natural isomorphisms above satisfy the coherence
relations of \cite{De} to define a genuine $B_{n+1}$ action on the
category $D^b(X)$.)
\end{theorem}

Thus we get braid group actions on derived categories of coherent
sheaves in the following cases.

\begin{itemize}
\item $A_n$-chains of $-2$-spheres $C_i$ (i.e. a sequence of $-2$-spheres
  with consecutive pairwise intersections a reduced point) in
  quasi-projective surfaces give actions of $B_{n+1}$ with
  $\F_i=\OO_{C_i}$. For instance the ALE spaces that are the
  resolutions of the $SL(2,\C)$ quotient singularities $\C^{\,2}/\mathbb
  Z_n$.
\item If, as in the previous example, we have an $A_n$-chain of
  $-2$-spheres $C_i$, we can twist instead by the corresponding line
  bundles $L_i=\OO(C_i)$. Simple exact sequences show that these also
  form an $A_n$-chain of spherical objects in $D^b(X)$ for $X$ a $K3$ surface.
\item Chains of surfaces $S_i$ in 3-folds, each with
  $\nu_S\cong\omega_S$ and $h^{0,1}(S)=0=h^{0,2}(S)$, give braid group
  actions (for $\F_i=\OO_{S_i}$)
  if $S_i\cap S_j=\emptyset$ for $|i-j|>1$ and $S_i\cap S_{i+1}$ is a
  $\Pee^1$-fibre of one surface and a $-2$-sphere in the
  other. Again such configurations arise in crepant resolutions of 3-fold
  singularities.
\item Taking $\F_i$ to be line bundles on an elliptic curve, with the
  degrees of $\F_i$ and $\F_j$ differing by $i-j$, we recover the
  original Fourier-Mukai transforms of \cite{Mu}, and from any two
  consecutive such line bundles we get an action of $B_3$,
  a central extension of $SL(2,\mathbb Z)$. The central
  generator $(T_1\,T_2)^6$ acts as the translation $[\,2\,]$ and the
  action is easily seen to be the $SL(2,\mathbb Z)$ action of
  Mukai. There are mirror, symplectic analogues of these relations for
  Dehn twists on tori in \cite{ST}.
\end{itemize}

\subsection*{Singularities}

Notice the chains of $-2$-spheres in surfaces are different in the
holomorphic and symplectic cases -- the former appear in resolutions
of singularities, the latter in smoothings of the same singularities. 
Are the corresponding braid group actions, on $D^b(X)$ and
$D^b(\mathrm{Fuk\,}(X))$ respectively, mirror dual\,?

In the compact case we know what the mirror to the Lagrangian spheres
should be in the presence of an SYZ fibration -- the Fourier-Mukai
transforms of their structure sheaves. If the spheres lie in the
elliptic fibres then their transforms are themselves, and the braid
group action on the symplectic side should be dual to the first
example listed above. If however the spheres are sections of the
fibration then their transforms are given by the line bundles
$\OO(C_i)$ (twisted by $\OO(-s)$, where $s$ is the image of the zero
section, but this need not concern us), and the correct dual is the
second example in the above list.

Either way mirror symmetry seems rather local in these cases, and
there are other cases where the mirror dual of the (symplectic)
smoothing of one singularity is the (holomorphic) resolution of
another. In fact this is the general proposal of \cite{Mo}: Morrison
suggests that moving towards the discriminant locus in the complex
structure moduli space of $X$ (i.e. degenerating $X$ to a singular
Calabi-Yau) should be mirror dual to moving to a ``boundary wall'' of
the complexified K\"ahler cone of $\widehat X$ (the annihilator of
a face of the Mori cone), thus inducing an extremal contraction of
$\widehat X$. Resolving the singularities of $X$ should
then be mirror dual to smoothing the contracted $\widehat X$. In
particular, in some generic situations, the smoothing of ordinary
double points (with their
Lagrangian $S^3$ vanishing cycles and corresponding Dehn twists on
$\widehat X$) should be mirror to small resolutions of other
ordinary double points (with exceptional $\Pee^1$ loci giving
corresponding twists $T_{\OO_{\Pee^1}}$ on $D^b(X)$).

Our results indicate that it might be reasonable to expect certain
additional properties of singularities that are dual in this
way. Namely, singularities whose smoothings have a Dynkin
diagrams of Lagrangian $S^3$ vanishing cycles should be dual to
singularities whose resolution has a similar diagram of (spherical
structure sheaves of) irreducible components of its exceptional set
(these form the nodes of the diagram, edges are provided by $\R\Hom$s).

For instance consider the smoothing of
the 3-fold singularity
$$
x^2+y^2+z^2+t^{2n}=0,
$$
which contains an $A_{2n-1}$-chain of Lagrangian $S^3$ vanishing
cycles of the singularity giving a braid group of (graded)
symplectomorphisms. Peturbing this into $n$ ordinary double points
$x^2+y^2+z^2+t^2=0$ (which we can do on the symplectic side, we are
only varying complex structure), analysing the
effect on homology and its mirror, and using Morrison's proposal, we
are led in \cite{ST} to ask whether the mirror should be given by
the following geometry, which we know leads to a braid group action.

\begin{prop}
Suppose we have a chain of Fano surfaces $\{S_{2k}\}_{k=1}^{n-1}$ in a
smooth 3-fold $X$, with the local Calabi-Yau condition that their
normal bundles are isomorphic to their canonical bundles. Suppose also
that each $S_{2k}$ contains two disjoint $(-1)$-$\Pee^1$s, $C_{2k-1}$ and
$C_{2k+1}$, in which it intersects $S_{2k-2}$ and $S_{2k+2}$ respectively
(and there are no more intersections, so only consecutive surfaces
intersect). Then the sheaves
$$
\F_{2k}=\OO_{S_{2k}}\quad\mathrm{and}\quad\F_{2k+1}=\OO_{C_{2k+1}}
$$
form an $A_{2n-1}$-chain, and so define an action of $B_{2n}$ on $D^b(X)$.
\end{prop}

In particular if we take the surfaces $S_{2k}$
to be $\Pee^2$s blown up in two points (giving the two exceptional curves
$C_{2k-1}$ and $C_{2k+1}$ which we think of as the mirrors of the
$(2k-1)$th and $(2k+1)$th vanishing cycles in the smoothing of
$x^2+y^2+z^2+t^{2n}=0$ according to Morrison's proposal) then there
is an extra $(-1)$-curve $C_{2k}$ in $S_{2k}$ -- the proper transform
of the line joining the two blow-up points -- which we can think of as
the mirror of the $2k$th vanishing cycle. Such a configuration is
easily shown \cite {ST} to arise in smooth
toric Calabi-Yau manifolds as the crepant resolution of a nasty
singularity that we would like to think of as the dual of
$x^2+y^2+z^2+t^{2n}=0$. 

Another relevant example is Arnold's strange duality
(see for instance \cite{Pi}), which is encompassed by mirror symmetry
for $K3$ surfaces according to the work of a number of people
(Aspinwall and Morrison, Kobayashi, Dolgachev, Ebeling, etc.).
To every isolated surface singularity on Arnold's list, described by
three numbers $b_1,\,b_2,\,b_3$, there is a natural $K3$
compactification $S$ of the singularity containing at infinity a chain of
$-2$-spheres with intersection configuration given by the following
Dynkin diagram $T(b_1,b_2,b_3)$:

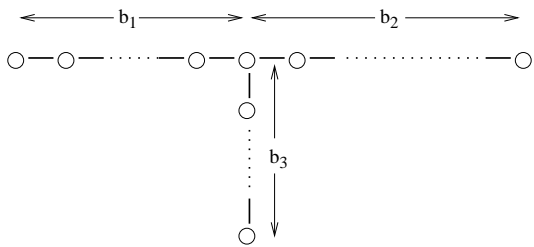
\begin{figure}[h]
\begin{center}
\includegraphics[angle=0,width=7cm]{dynk.eps}
\end{center}
\end{figure}
Here the circles represent $-2$-spheres, and edges give intersections
of the spheres of intersection number 1, and the central $-2$-sphere
is counted in each $b_i$.
The corresponding intersection matrix is not negative definite for
the numbers $b_i$ in the list, so the spheres cannot be completely
contracted, though they can be contracted to a smooth $\Pee^1$ with three
points on it that are singular points of the surface (this surface is
Pinkham's original compactification of the affine surface
singularity).

Let $\{c_i\}_{i=1}^3$ be the numbers dual to the $b_i$\,s in strange
duality (i.e. these are the $b_i$\,s associated to the dual
singularity). Then the intersection matrix of the smoothing of the
original singularity has
vanishing cycles given by $T(c_i)\oplus H$ (where $H$ is the
hyperbolic $\left(\!\!\!\begin{array}{cc} 0 & 1 \\ 1 & 0
\end{array}\!\!\!\right)$). Together the $-2$-spheres plus the vanishing
cycles give all of the homology of the smoothed $K3$ surface:
$H_2(K3)\cong T(b_i)\oplus T(c_i)\oplus H$. The duality swaps the
$b_i$\,s and $c_i$\,s, taking the homology $T(c_i)\oplus H$
generated by the vanishing cycles to the holomorphic homology
$T(c_i)\oplus H$ in the mirror $K3$ (given by the resolution cycles
and the hyperbolic $H_0\oplus H_4$).

Thus we can think of mirror symmetry as replacing the smoothing of one
singularity by the resolution of another (though it is not an isolated
singularity, it is the $\Pee^1$ with 3 surface singularities on
it). Here we are not thinking of the duality as merely swapping
Arnold's singularities, but as a rather more global phenomenon on $K3$
surfaces, which therefore has to include the chains of $-2$-spheres at
infinity.

There is a generalised braid group action on the derived category of
the $K3$ surface given by the $T(b_1,b_2,b_3)$ $-2$-sphere configuration.
Dually, there is probably the same group of Dehn twists around the
Lagrangian vanishing cycles of the mirror $K3$, though it seems not to
be known whether they can be found in the \emph{geometric} intersection
configuration $T(b_1,b_2,b_3)$ (they may have many more geometric
intersections than their topological intersections suggest).

\section{Mutations} \label{mut}

The formula (\ref{twist}) for our twist is familiar to
algebraists, in tilting theory, and those who work on exceptional
bundles on Fano manifolds -- see e.g. \cite{Ru}. There the twists are
called mutations, and act on certain modules over algebras (similar to
those described in the next section) built from the bundles
(they cannot give equivalences of derived categories by the result of
\cite{BO}). The bundles $\F$ that one twists are also those with
minimal $\Ext^*(\F,\F)$; in the case of Fano manifolds this means $\F$
is simple and has no higher Exts at all, and is called
\emph{exceptional}.

There are braid group actions of such twists on exceptional collections of
bundles, but the relation with our work is far from clear, and it is
possible there is none. Here we shall simply note a relationship
between exceptional objects on Fano manifolds and spherical objects on
Calabi-Yaus, motivated by two examples of \cite{Ku}.

\begin{definition} We say that a map $f:\,X\to Y$ from a Calabi-Yau
$N$-fold $X$ to a smooth projective Fano $M$-fold $Y$, of codimension
$c=N-M$ (which may be of any sign or zero), is \emph{simple} if
$\R f_*\OO_X$ is made up from $\OO_Y$ and $\omega_Y$ in the sense that
\begin{itemize}
\item{$c>0\qquad$} $R^if_*\OO_X=\left\{ \begin{array}{cl}
\OO_Y & i=0 \\ 0 & i\ne0,\,c \\ \omega_Y & i=c \end{array} \right.$
\item{$c=0$\qquad} $R^if_*\OO_X=\left\{ \begin{array}{cl}
\OO_Y\oplus\omega_Y & i=0 \\ 0 & i\ne0 \end{array} \right.$
\item{$c<0$\qquad} There is an exact triangle $\omega_Y\,[-c]\to\R
f_*\OO_X\to\OO_Y$.
\end{itemize}
\end{definition}   

Simple examples of maps $f$ are often of this type. For instance, for
Calabi-Yaus fibred over a Fano base with generic fibre $F$ such that
$h^{0,i}(F)=0,\ 0<i<c$, relative Serre duality shows that
the projection is simple in this sense.

Examples with $c=0$ are given by Calabi-Yaus double covering Fanos,
branched over a double anticanonical divisor, while for $c=-1$ we have
a Calabi-Yau anticanonical divisor in a Fano manifold.

\begin{theorem}
Suppose $f:\,X\to Y$ is a simple map, as defined above, and
$\F\in D^b(Y)$ is exceptional $(\Ext^i(\F,\F)=\C$ for $i=0$,
and $0$ for $i\ne0)$. Then $\Lf\F\in D^b(X)$ is spherical.
\end{theorem}

In the other direction, i.e. maps from a Fano to a Calabi-Yau,
something can only be said in the case of a Fano divisor in a
(locally) Calabi-Yau manifold. Namely,

\begin{theorem}
Suppose that $\iota:\,Y\subset X$ is a smooth Fano divisor with normal
bundle $\nu_Y=\omega_Y$ in a quasi-projective scheme $X$. If
$\F\in D^b(Y)$ is exceptional then $\iota_*\F\in D^b(X)$ is spherical.
\end{theorem}

This takes care of most of the examples of spherical sheaves given
until now (by taking the exceptional sheaf to be $\OO_Y$), and
provides many more by pushing forward the exceptional collections of
\cite{Ru} to (local) Calabi-Yau manifolds.

\section{Fidelity} \label{fid}

Finally we briefly mention what is in many ways the main result of
\cite{KS,ST}, namely that in dimension $N\ge2$, the $B_{n+1}$
actions given by $A_n$-chains of spherical objects are
\emph{faithful}.

To do this it is clearly enough to show the induced $B_{n+1}$ action on the
differential graded modules $\bigoplus_i\R\Hom(\F_i,\G)$, in the
derived category of the
differential graded algebra $\bigoplus_{ij}\R\Hom(\F_i,\F_j)$, is
faithful. (Here the algebra and module structures are the obvious
ones; replacing each $\F_i$ by a finite resolution of injectives,
which one can prove is possible, it is just composition of morphisms.)

In fact a difficult result of \cite{KS} is that the induced action on
homology, i.e. the action on the graded modules
$\bigoplus_{ik}\Ext^k(\F_i,\G)$ in the derived category of the graded algebra
$\bigoplus_{ijk}\Ext^k(\F_i,\F_j)$, is faithful. (Since the $\F_i$\,s
form an $A_n$-chain these modules and algebras take a standard form,
and the braid group action is the one considered in \cite{KS}.)
However this is not enough to prove faithfulness of the action at the
level of differential graded modules. The following result, though,
is sufficient to provide a proof.

\begin{theorem}
For $\{\F_i\}_{i=1}^n$ an $A_n$-chain of $N$-spherical objects,
$N\ge2$, the graded algebra $A=\bigoplus_{ik}\Ext^k(\F_i,\F_i)$ is
\emph{intrinsically formal}. That is, any differential graded algebra
with $A$ as its cohomology is quasi-isomorphic to $A$.
\end{theorem}

This is proved by a lot of non-commutative obstruction theory that I will
not go into here. That such machinery is really necessary is seen in
the following example of non-faithfulness in dimension $N=1$.

Let $T$ be an elliptic curve, $L$ be a non-trivial degree
zero line bundle, and denote by $\OO_p$ the structure sheaf of a
point $p\in T$. Notice that $L$ induces an automorphism $\phi_L$
of $T$: identifying $T$ with its degree one line bundles, $\phi_L$ is
given by tensoring with $L$.

The sheaves $\OO=\OO_T,\ \OO_p$ and $L$ form an $A_3$-chain of
$1$-spherical sheaves.

\begin{theorem}
The action of $T'_L T_\OO$ is the action induced by the automorphism
$\phi_L$ by pullback. In particular if $L$ has order two
$(L^2=\OO)$ then $(T_L^{-1}T_\OO)^2\cong\id.$
\end{theorem}

\bibliographystyle{amsalpha}

\begin{thebibliography}{SYZ}


\bibitem[AP]{AP} D. Arinkin and A. Polishchuk,
\textit{Fukaya category and Fourier transform}, preprint
math/9811023.

\bibitem[BO]{BO} A.\,I. Bondal and D.\,O. Orlov,
\textit{Semiorthogonal decomposition for algebraic varieties}, preprint
alg-geom/9506012.

\bibitem[Br]{Br} T. Bridgeland,
\textit{Equivalences of triangulated categories and Fourier-Mukai transforms},
Bull. Lond. Math. Soc. \textbf{31} (1999), 25--34.

\bibitem[De]{De} P. Deligne,
\textit{Action du groupe des tresses sur une cat\'egorie},
Invent. Math. \textbf{128} (1997), 159--175.

\bibitem[KS]{KS} M. Khovanov and P. Seidel,
\textit{Quivers, Floer cohomology, and braid group actions},
in preparation, 1999.

\bibitem[K]{K} M. Kontsevich, 
\textit{Homological Algebra of Mirror Symmetry}, International Congress
of Mathematicians, Z\"urich 1994. Birkh\"auser, 1995.    

\bibitem[Ku]{Ku} S.\,A. Kuleshov,
\textit{Exceptional bundles on $K3$ surfaces}, in \cite{Ru}.

\bibitem[Ma]{Ma} A. Maciocia,
\textit{Generalized Fourier-Mukai transforms}, J. Reine Angew. Math.
\textbf{480} (1996), 197--211.

\bibitem[Mo]{Mo} D.\,R. Morrison,
\textit{Through the looking glass}, Mirror Symmetry III (Montreal 1995),
AMS/IP Stud. Adv. Math. 10 (1999), 263--277.

\bibitem[Mu]{Mu} S. Mukai,
\textit{Duality between $D(X)$ and $D(\widehat X)$ with its application to
Picard sheaves}, Nagoya Math. Jour. \textbf{81} (1981), 153--175.

\bibitem[Or]{Or} D.\,O. Orlov,
\textit{Equivalences of derived categories and K3 surfaces}, preprint
alg-geom/9606006.

\bibitem[Pi]{Pi} H. Pinkham,
\textit{Singularit\'es exceptionnelles, la dualit\'e \'etrange d'Arnold,
et les surfaces $K3$}, C. R. Acad. Sc. Paris \textbf{284\,A} (1977)
615--618.

\bibitem[Ru]{Ru} A.\,N. Rudakov et al.,
\textit{Helices and vector bundles: Seminaire Rudakov}, LMS Lecture Note
Series 148, Cambridge University Press, 1990.

\bibitem[S]{S} P.\,S. Seidel,
\textit{Lagrangian two-spheres can be symplectically knotted}, preprint
math/9803083.

\bibitem[ST]{ST} P.\,S. Seidel and R.\,P. Thomas,
\textit{Braid group actions on derived categories of sheaves},
pre\-print 1999.

\bibitem[SYZ]{SYZ} A. Strominger, S.-T. Yau and E. Zaslow,
\textit{Mirror Symmetry is T-Duality}, Nucl. Phys. \textbf{B479}
(1996), 243--259.

\bibitem[Ty]{Ty} A.\,N. Tyurin,
\textit{Geometric quantization and mirror symmetry},
preprint math.AG/9902027.

\bibitem[Va]{Va} C. Vafa,
\textit{Extending Mirror Conjecture to Calabi-Yau with Bundles},
preprint hep-th/ 9804131.

\end{thebibliography}

\end{document}